\documentstyle[11pt, a4wide, leqno]{article}

\def\bea{\begin{eqnarray}}
\def\eea{\end{eqnarray}}
\def\beann{\begin{eqnarray*}}
\def\eeann{\end{eqnarray*}}
\def\beq{\begin{equation}}
\def\eeq{\end{equation}}
\def\ba{\begin{array}}
\def\ea{\end{array}}
\def\ben{\begin{enumerate}}
\def\een{\end{enumerate}}

\newtheorem{th}{Theorem}[section]

\newtheorem{lem}[th]{Lemma}

\newtheorem{pro}[th]{Proposition}

\newtheorem{co}[th]{Corollary}

\def\e  {\epsilon}

\date{}

\begin{document}
%\rightline{}

\thispagestyle{empty}
\begin{center}
\vspace {.7cm} {\Large {\bf Connections with torsion, parallel
spinors\\ and geometry of
 Spin(7) manifolds} }
 \vskip 1.0truecm
{\large{\bf
Stefan Ivanov}} \vskip 1.0truecm

e-mail: ivanovsp@fmi.uni-sofia.bg
%{\normalsize{\sl  Department of Mathematics}}
%\\ {\normalsize{\sl University of Sofia}}
%\\ {\normalsize{\sl \lq\lq St. Kl. Ohridski''}}
\end{center}
\begin{abstract}
We show that on every $Spin(7)$-manifold there always
exists a unique linear connection with totally skew-symmetric torsion preserving
 a nontrivial spinor and the $Spin(7)$
structure. We express its torsion and the Riemannian scalar
curvature in terms of the fundamental 4-form. We present an
explicit formula for the Riemannian covariant derivative of the
fundamental 4-form in terms of its exterior differential. We show
the vanishing of the $\hat A$-genus and obtain a linear relation
between Betti numbers of a compact $Spin(7)$ manifold which is
locally but not globally conformally equivalent to a space with
closed fundamental 4-form. A general solution to the Killing
spinor equations is presented.
\\[15mm]
{\bf Running title:} Geometry of $Spin(7)$-manifolds
\\[2mm]
{\bf Subj. Class.:} Special Riemannian manifolds, Spin geometry, String theory
\\[2mm]
${\bf MS}$ {\bf classification: } 53C25; 53C27; 53C55; 81T30
\\[2mm]
\end{abstract}

%\newpage
\section{Introduction}

Riemannian manifolds admitting parallel spinors with respect to a
metric connection with totally skew-symmetric torsion recently
become a subject of interest in theoretical and mathematical
physics. One of the main reasons is that the number of preserving
supersymmetries in string theory depends essentially on the number
of parallel spinors. In 10-dimensional string theory, the Killing
spinor equations in the string frame can be written in the
following way \cite{Stro}, (see eg \cite{IP1,IP,FI})
\beq\label{ks1} \nabla \psi=0, \eeq \beq\label{ks2} (d\Psi -
\frac{1}{2}H)\cdot\psi=0, \eeq where $\Psi$ is a scalar function
called the dilation, $H$ the 3-form field strength, $\psi$ a
spinor field and $\nabla$ a metric connection with totally
skew-symmetric torsion $T=H$. The number of preserving
supersymmetries is determined by the number of  solutions of these
equations.

The existence of a  parallel spinor imposes restrictions on the
holonomy group since the spinor holonomy representation has to
have a fixed point. In the case of torsion-free metric connections
(Levi-Civita connections) the possible Riemannian holonomy groups
are known to be SU(n), Sp(n), $G_2$, Spin(7) \cite{Hit,Wang}. The
Riemannian holonomy condition imposes strong restrictions on the
geometry and leads to considerations of Calabi-Yau manifolds,
hyper-K\"ahler manifolds, parallel $G_2$-manifolds, parallel
$Spin(7)$ manifolds. All of them are of great interest in
mathematics (see \cite{J2} for precise discussions) as well as in
high-energy physics, string theory \cite{Pol}.

It just happens that the geometry of these spaces is too
restrictive for various questions in string theory
\cite{MN,ST,GKMW}. It seems that a 'nice' mathematical
generalization of Calabi-Yau manifolds, hyper-K\"ahler manifolds,
parallel $G_2$-manifolds, parallel $Spin(7)$ manifolds is to
consider linear connections with skew-symmetric torsion and
holonomy contained in $SU(n), Sp(n), G_2, Spin(7)$.

A remarkable fact is that the existence (in small dimensions) of a
parallel spinor with respect to a metric connection with
skew-symmetric torsion determines the connection in a unique way
in the cases where its holonomy group is a subgroup of $SU, Sp,
G_2$ provided additional differential conditions on the structure
are fulfilled  \cite{Stro,FI}. The uniqueness property leads to
the idea that it is worth to study the geometry of such a
connection with torsion, besides its interest in physics
\cite{Pol,GKMW}, for purely mathematical reasons expecting to get
information about the curvature of the metric, Betti numbers,
Hodge numbers, $\hat A$-genus, etc. In fact, a connection with
skew symmetric torsion preserving a given complex structure on a
Hermitian manifold was used by Bismut \cite{Bis} to prove a local
index formula for the Dolbeault operator when the manifold is not
K\"ahler. Following this idea, a vanishing theorem for the
Dolbeault cohomology on a compact Hermitian non-K\"ahler manifold
was found \cite{AI,IP,IP1}.

In this paper we study the existence of parallel spinors with
respect to a metric connection with skew-symmetric torsion in
dimension 8 (for dimensions 4,5,6,7 see \cite{Stro,DI,IP,FI,FI1}).
The first consequence is that the manifold should be a $Spin(7)$
manifold, i.e. its structure group can be reduced to the group
$Spin(7)$. This is because the Euler characteristic ${\cal
X}(S_{\pm})$ of at least one of the (negative $S_-$ or positive
$S_+$) spinor bundles vanishes and therefore the structure group
can be reduced to $Spin(7)$ \cite{LM}. Surprisingly, we discover
that the converse is always true in dimension 8. We show that the
existence of a connection with totally skew-symmetric torsion
preserving a spinor in dimension 8 is unobstructed, i.e. on every
$Spin(7)$ 8-manifold there always exists a unique linear
connection with totally skew-symmetric torsion preserving a
nontrivial spinor ie with holonomy contained in $Spin(7)$. This
phenomena does not occur in the cases of holonomy groups $SU, Sp,
G_2$ (see the end of the paper). We find a formula for the torsion
3-form and for the Riemannian scalar curvature in terms of the
fundamental 4-form. Our main result is the following
\begin{th}\label{th2}
Let $(M,g,\Phi)$ be an 8-dimensional $Spin(7)$ manifold with
fundamental 4-form $\Phi$.
\item
i). There always exists a unique linear connection $\nabla$
preserving the $Spin(7)$ structure, $\nabla\Phi=\nabla g=0$, with
totally skew-symmetric torsion $T$ given by \beq\label{n4}
T=-\delta \Phi - \frac{7}{6}*(\theta\wedge\Phi),\quad \theta =
\frac{1}{7}*(\delta\Phi\wedge\Phi). \eeq On any $Spin(7)$ manifold
there exists a $\nabla$-parallel spinor $\phi$ corresponding to
the fundamental form $\Phi$ and the Clifford action of the torsion
3-form on it is \beq\label{b2} T\cdot\phi= - \frac{7}{6}\theta
\cdot\phi. \eeq
\item
ii). The Riemannian scalar curvature $Scal^g$ and the scalar
curvature $Scal$ of the $Spin(7)$ connection $\nabla$ are given in
terms of the fundamental 4-form $\Phi$ by
\begin{eqnarray}\label{sc1}
Scal^g = \frac{49}{18}||\theta||^2 -\frac{1}{12}||T||^2+
\frac{7}{2}\delta\theta, \qquad \qquad Scal =
\frac{49}{18}||\theta||^2
-\frac{1}{3}||T||^2+\frac{7}{2}\delta\theta.
\end{eqnarray}
\end{th}
The proof relies on our explicit formula expressing the covariant
derivative of the fundamental 4-form $\Phi$ with respect to the
Levi-Civita connection in terms of the exterior derivative of
$\Phi$. The existence of such a relation was discovered by
R.L.Bryant \cite{Br} in his proof that the holonomy group of the
Levi-Civita connection is contained in $Spin(7)$ iff $d\Phi=0$
(see also \cite{Sal}). We prove ii) using the
Schr\"odinger-Lichnerowicz formula for the connection with torsion
established in \cite{FI} and the special properties of the
Clifford action on the special spinor $\phi$.

In the compact case, we use  the formula for the Riemannian scalar
curvature to show that the Yamabe constant of one of the two
classes of $Spin(7)$ manifolds according to Fernandez
classification \cite{F} is strictly positive. Applying the
Atiyah-Singer index theorem \cite{AS}, as well as the Lichnerowicz
vanishing theorem \cite{Li}, we find a linear relation between the
Betti numbers and show that the Euler characteristic is equal to 3
times the signature.

In the last section we give necessary and sufficient conditions
for the existence of a solution to the Killing spinor equations
(\ref{ks1}), (\ref{ks2}) in an 8-dimensional manifold. We apply
our general formula for the torsion of the connection admitting a
parallel spinor to the second Killing spinor equation. As a
consequence, we obtain a formula for the field strength (torsion)
of a solution to both Killing spinor equations in terms of the
fundamental 4-form. We discover a relation between a solution to
both Killing spinor equations with non-constant dilation and the
conformal transformations of the $Spin(7)$ structures. In fact we
show that the dilation function arises geometrically (from the Lee
form of the structure) and can be interpreted as a conformal
factor. Our analysis on the two Killing spinor equations in
dimension 8 shows that the physics data (field strength $H$ and
the dilation function $\Psi$) are determined completely by the
properties of the parallel spinor or equivalently by the geometry
of the corresponding fundamental 4-form.

{\bf Acknowledgements} The research was done during the author's
visit at the Abdus Salam International Centre for Theoretical
Physics, Trieste, Italy, Fall 2001. The author thanks the Abdus
Salam ICTP for its support and the excellent facilities. The
author is a member of the EDGE, Research Training Network
HPRN-CT-2000-00101, supported by the European Human Potential
Programme. The research is partially supported by Contract MM
809/1998 with the Ministry of Science and Education of Bulgaria,
Contract 586/2002 with the University of Sofia "St. Kl. Ohridski".

\section{General properties of $Spin(7)$ manifold}
We recall some notions of $Spin(7)$ geometry.

Let us consider ${\cal R}^8$ endowed with an orientation and its
standard inner product $<,>$. Let $\{e_0,...,e_7\}$ be an oriented
orthonormal basis. We shall use the same notation for the dual
basis. We denote by $e_{ijkl}$ the monomial $e_i\wedge$$ \e_j
$$\wedge$$ \e_k$$ \wedge$$ \e_l$. Consider the 4-form $\Phi$ on
${\cal R}^8$ given by
\begin{eqnarray}\label{1}
\Phi &=&e_{0123} + e_{0145} + e_{0167}+ e_{0246} - e_{0257} -
e_{0347} - e_{0356} \\ \nonumber &+& e_{4567} +  e_{2367}+
e_{2345} + e_{1357}- e_{1346}- e_{1247} - e_{1256} \quad
.\nonumber
\end{eqnarray}
The 4-form  $\Phi$ is self-dual $*\Phi=\Phi$, where $*$ is the
Hodge $*$-operator and the 8-form $\Phi\wedge\Phi$ coincides with
the volume form of ${\cal R}^8$. The subgroup of $GL(8,R)$ which
fixes $\Phi$ is isomorphic to the double covering $Spin(7)$ of
$SO(7)$ \cite{HL}. Moreover, $Spin(7)$ is a compact
simply-connected Lie group of dimension 21 \cite{Br}. The 4-form
$\Phi$ corresponds to a real spinor $\phi$ and therefore,
$Spin(7)$ can be identified as the isotropy group of a non-trivial
real spinor.

A 3-fold vector cross product $P$ on ${\cal R}^8$ can be defined
by $<P(x\wedge y\wedge z),t>= \Phi(x,y,z,t)$, for $x,y,z,t \in
{\cal R}^8$. Then $Spin(7)$ is also characterized by $$
Spin(7)=\{a\in O(8)| P(ax\wedge ay\wedge az)=P(x\wedge y\wedge z),
x,y,z \in {\cal R}^8\}.$$ The inner product $<,>$ on ${\cal R}^8$
can be reconstructed from $\Phi$ \cite{F,Gr}, which corresponds
with the fact that $Spin(7)$ is a subgroup of $SO(8)$.

A {\it $Spin(7)$ structure} on an 8-manifold $M$ is by definition
a reduction of the structure group of the tangent bundle to
$Spin(7)$; we shall also say that $M$ is a {\it $Spin(7)$
manifold}. This can be described geometrically by saying that
there is a 3-fold vector cross product $P$ defined on $M$, or
equivalently there exists a nowhere vanishing differential 4-form
$\Phi$ on $M$ which can be locally written as (\ref{1}). The
4-form $\Phi$ is called the {\it fundamental form} of the
$Spin(7)$ manifold $M$ \cite{Bo}.

Let $(M,g,\Phi)$ be a $Spin(7)$ manifold. The action of $Spin(7)$
on the tangent space gives an action of $Spin(7)$ on
$\Lambda^k(M)$ and so the exterior algebra splits orthogonally
into components, where $\Lambda^k_l$ corresponds to an irreducible
representation of $Spin(7)$ of dimension $l$ \cite{F,Br}:
$$\Lambda^1(M)=\Lambda^1_8, \quad \Lambda^2(M) = \Lambda^2_7\oplus
\Lambda^2_{21}, \quad
\Lambda^3(M)=\Lambda^3_8\oplus\Lambda^3_{48}, $$ $$
\Lambda^4(M)=\Lambda^4_+(M)\oplus \Lambda^4_-(M), \quad
\Lambda^4_+(M)=\Lambda^4_1\oplus\Lambda^4_7\oplus\Lambda^4_{27},
\quad \Lambda^4_-=\Lambda^4_{35}; $$
where $\Lambda^4_{\pm}(M)$
are the $\pm$-eigenspaces of $*$ on $\Lambda^4(M)$ and
\cite{Bo,Br,Sal} $$ \Lambda^2_7 = \{\alpha \in \Lambda^2(M) |
*(\alpha\wedge\Phi)=3\alpha\}, \quad \Lambda^2_{21} = \{\alpha \in
\Lambda^2(M)|*(\alpha\wedge\Phi)=-\alpha\} $$ $$ \Lambda^3_8 =
\{*(\beta\wedge\Phi) | \beta \in \Lambda^1(M)\}, \quad
\Lambda^3_{48} = \{\gamma \in \Lambda^3(M) | \gamma\wedge\Phi=0\},
\quad \Lambda^4_1 = \{f\Phi | f\in {\cal F(M)}\} $$ The Hodge star
$*$ gives an isometry between $\Lambda^k_l$ and $\Lambda^{8-k}_l$.

If $(M,g,\Phi)$ is a $Spin(7)$ manifold, then $M$ is orientable
and spin, with preferred spin structure and orientation. If
$S=S_+\oplus S_-$ is the spin bundle of $M$, then there are
natural isomorphisms $S_+\equiv \Lambda^0_1\oplus\Lambda^2_7$ and
$S_-\equiv \Lambda^1_8$ (see eg \cite{J2}).

In general, not every 8-dimensional Riemannian spin manifold $M^8$
admits a $Spin(7)$ structure. We explain the precise condition
\cite{LM}. Denote by $p_1(M), p_2(M), {\cal X}(M), {\cal
X}(S_{\pm})$ the first and the second Pontrjagin classes, the
Euler characteristic of $M$ and the Euler characteristic of the
positive and the negative spinor bundles, respectively. It is well
known \cite{LM} that a spin 8-manifold admits a $Spin(7)$
structure if and only if ${\cal X}(S_+)=0$ or ${\cal X}(S_-)=0$.
The latter conditions are equivalent to \cite{LM} \beq\label{c1}
p_1^2(M)-4p_2(M)+ 8{\cal X}(M)=0, \eeq for an appropriate choice
of the orientation.

Let us recall that a $Spin(7)$ manifold $(M,g,\Phi)$ is said to be
parallel (torsion-free \cite{J2}) if the holonomy of the metric
$Hol(g)$ is a subgroup of $Spin(7)$. This is equivalent to saying
that the fundamental form $\Phi$ is parallel with respect to the
Levi-Civita connection $\nabla^g$ of the metric $g$. Moreover,
$Hol(g)\subset Spin(7)$ if and only if $d\Phi=0$ \cite{Br} (see
also \cite{Sal}) and any parallel $Spin(7)$ manifold is Ricci flat
\cite{Bo}.

According to the Fernandez classification \cite{F}, there are
4-classes of $Spin(7)$ manifolds  obtained as irreducible
representations of $Spin(7)$ of the space $\nabla^g\Phi$.
Following \cite{C1} we consider the 1-form $\theta$ defined by
\beq\label{c2} 7\theta = -*(*d\Phi\wedge\Phi)=*(\delta\Phi\wedge
\Phi) \eeq We shall call the 1-form $\theta$ {\it the Lee form} of
a given $Spin(7)$ structure.

The 4 classes of Spin(7) manifolds in the Fernandez classification
can be described in terms of the Lee form as follows \cite{C1}:
$W_0 : d\Phi=0; \quad W_1 : \theta =0; \quad W_2 : d\Phi =
\theta\wedge\Phi; \quad W : W=W_1\oplus W_2.$

We shall call a $Spin(7)$ structure of the class $W_1$ (ie
$Spin(7)$ structures with zero Lee form)  {\it a balanced}
$Spin(7)$ structure.

In \cite{C1} Cabrera shows that the Lee form of a $Spin(7)$
structure in the class $W_2$ is closed and therefore such a
manifold is locally conformally equivalent to a parallel $Spin(7)$
manifold  and it is called {\it locally conformally parallel}. If
the Lee form is not exact (i.e. the structure is not globally
conformally parallel), we shall call it {\it strict locally
conformally parallel}. We shall see later (section 8) that these
spaces have very different topology than parallel ones.

Coeffective cohomology and coeffective numbers of Riemannian
manifolds with $Spin(7)$ structure are studied in \cite{Ug}.

\section{Examples:}

Examples of $Spin(7)$ manifolds are constructed relatively recently.

The first known explicit example of complete parallel $Spin(7)$
manifold with $Hol(g)=Spin(7)$ was constructed by Bryant and
Salamon \cite{BS,Gibb} on the total space of the spin bundle over
the 4-sphere.

The first compact examples of parallel $Spin(7)$ manifolds with
$Hol(g)=Spin(7)$ were constructed by Joyce\cite{J1,J2} by
resolving the singularities of the orbifold $T^8/\Gamma$ for
certain discrete groups $\Gamma$.

Most examples of $Spin(7)$ manifolds in the Fernandez
classification are constructed by using certain $G_2$-manifolds.
We recall that a $G_2$-manifold $N$ is a 7-dimensional manifold
whose structure group can be reduced to the exceptional group
$G_2$ or equivalently, there exists on $N$ a distinguished
associative 3-form $\gamma$. A $G_2$-manifold is said to be {\it
nearly parallel, cocalibrated of pure type, calibrated} if
$d\gamma=const.*\gamma; \quad \delta\gamma=0,
d\gamma\wedge\gamma=0; d\gamma=0$, respectively \cite{FG}.

Any 8-manifold of type $M=S^1\times N$ possesses a $Spin(7)$
structure defined by \cite{Sal,Ug,Ca} $\Phi=\eta\wedge\gamma +
*\gamma$, where $\eta$ is a non-zero 1-form on $S^1$. The induced
$Spin(7)$ structure on $M$ is \cite{Ca}

i) a strict locally conformally parallel if the
$G_2$ structure is nearly parallel;

ii) a balanced one if the $G_2$-structure is cocalibrated of pure
type or calibrated or belongs to the direct sum of these classes.

There are many known examples of compact nearly parallel
$G_2$-manifolds: $S^7$ \cite{FG}, $SO(5)/SO(3)$ \cite{BS,Sal}, the
Aloff-Wallach spaces $N(g,l)=SU(3)/U(1)_{g,l}$ \cite{CMS} any
Einstein-Sasakian and any 3-Sasakian space in dimension 7
\cite{FK,FKMS}, some examples coming from 7-dimensional 3-Sasaki
manifolds \cite{FKMS,GS}, the 3-Sasakian non-regular spaces
$S(p_1,p_2,p_3)$ \cite{6,7}, compact nearly parallel
$G_2$-manifolds with large symmetry groups are classified recently
in \cite{FKMS}. The product of each of these spaces by $S^1$ gives
examples of strict locally conformally parallel $Spin(7)$
structures.

Any minimal hypersurface $N$ in $R^8$ possesses  a cocalibrated
structure of pure type $G_2$ \cite{FG} and therefore $M=N\times
S^1$ has a balanced $Spin(7)$ structure described above.

More general, any principle fibre bundle with one dimensional
fibre over a $G_2$-manifold carries a $Spin(7)$ structure
\cite{C1}. In this way, a balanced $Spin(7)$ structure arises on a
principle circle bundle over a 7-dimensional torus $T^7$
considered as a $G_2$-manifold \cite{C1}.

\section{Conformal transformations of $Spin(7)$ structures}

We need the next result which is essentially established in \cite{F}.
\begin{pro}\label{pro1}\cite{F}
Let $\bar{g}=e^{2f}g, \quad \bar{\Phi}=e^{4f}\Phi$  be a conformal change of the
given $Spin(7)$ structure $(g,\Phi)$ and $\bar{\theta}, \quad \theta$
are the corresponding Lee 1-forms, respectively. Then
\beq\label{2}
\bar{\theta}= \theta + 4df
\eeq
\end{pro}
{\bf Proof.} We have \cite{F} $vol._{\bar{g}}=e^{8f}vol._g, \quad
d\bar{\Phi} = e^{4f}(4df\wedge\Phi + d\Phi)$. We calculate
$$\bar{*}d\bar{\Phi}=e^{4f}(*d\Phi + 4*(df\wedge\Phi)), \quad
\bar{*}d\bar{\Phi}\wedge\bar{\Phi}=e^{8f}(*d\Phi\wedge\Phi +
28*df),$$ where we used the identity $
*(\Phi\wedge\gamma)\wedge\Phi=7*\gamma, \quad \gamma \in
\Lambda^1(M)$. We obtain consequently that
$\bar{\theta}=-\frac{1}{7}
\bar{*}(\bar{*}d\bar{\Phi}\wedge\bar{\Phi})=-\frac{1}{7}
*(*d\Phi\wedge\Phi) -4 *^2df= \theta +4df.$ \hfill {\bf Q.E.D.}

More generally, we have
\begin{co}\label{co1}
If the Lee 1-form is closed, then the $Spin(7)$ structure is
locally conformal to a balanced $Spin(7)$ structure.
\end{co}
Proposition~\ref{pro1} allows us to find a distinguished $Spin(7)$ structure on a
 compact 8-dimensional $Spin(7)$ manifold.
\begin{th}\label{thg}
Let $(M^8,g,\Phi)$ be a compact 8-dimensional $Spin(7)$ manifold.
Then there exists a unique (up to homothety) conformal $Spin(7)$
structure $g_0=e^{2f}g, \Phi_0=e^{4f}\Phi$ such that the
corresponding Lee 1-form is coclosed, $\delta_0\theta_0=0$.
\end{th}
{\bf Proof.}  We shall use  the Gauduchon theorem for the
existence of a distinguished metric (Gauduchon metric) on a
compact Hermitian or Weyl manifold \cite{G1,G2}. We shall use the
expression of this theorem in terms of a Weyl structure (see
\cite{tod}, Appendix 1). We consider the Weyl manifold
$(M^8,g,\theta,\nabla^W)$ with the Weyl 1-form $\theta$ where
$\nabla^W$ is a torsion-free linear connection on $M^8$ determined
by the condition $\nabla^Wg=\theta\otimes g$. Applying the
Gauduchon theorem we can find in a unique way a conformal metric
$g_0$ such that the corresponding Weyl 1-form is coclosed with
respect to $g_0$. The key point is that by Proposition~\ref{pro1}
the Lee 1-form transforms under conformal rescaling according to
(\ref{2}) which is exactly the transformation of the Weyl 1-form
under conformal rescaling of the metric $\bar g=e^{4f}g$. Thus,
there exists a unique (up to homothety) conformal $Spin(7)$
structure $(g_0,\Phi_0)$ with coclosed Lee 1-form. \hfill {\bf
Q.E.D.}

We shall call the $Spin(7)$ structure with coclosed Lee 1-form
{\it the Gauduchon $Spin(7)$ structure}.

\begin{co}
Let $(M,g,\Phi)$ be a compact $Spin(7)$ manifold and $(g,\Phi)$ be
the Gauduchon structure. Then the following formula holds
\hspace{2mm}
$
*\left(d\delta\Phi\wedge\Phi\right)=-||d\Phi||^2.
$
\end{co}
{\bf Proof.}
Using (\ref{c2}), we calculate that
$
0=7\delta\theta=*d(*d\Phi\wedge\Phi)=*\left(d\delta\Phi\wedge\Phi
-*d\Phi\wedge d\Phi\right)=\\ =*\left(d\delta\Phi\wedge\Phi
+||d\Phi||^2.vol\right).$ \hfill {\bf Q.E.D.}

\begin{co}\label{cof}
On a compact $Spin(7)$ manifold with closed Lee form the first
Betti number $b_1\ge 1$ provided the Gauduchon $Spin(7)$ structure
is not balanced. In particular, on any strict locally conformally
parallel $Spin(7)$ manifold, $b_1\ge 1$.
\end{co}

\section{A formula for the covariant derivative of the fundamental form}

In \cite{Br} R.L. Bryant proved that on a $Spin(7)$ manifold
$(M,g,\Phi)$ the holonomy group $Hol(g)$ of the metric $g$ is
contained in $Spin(7)$ iff the fundamental form $\Phi$ is closed
i.e. $\nabla^g\Phi=0$ is equivalent to $d\Phi=0$. This shows that
there is an identification of $\nabla^g\Phi$ and $d\Phi$ (see also
\cite{Sal}). The aim of this section is to give an explicit
formula.

Let $\gamma$ be an 1-form, $\gamma\in\Lambda^1(M)$. We use the
same notation for the dual vector field via the metric and denote
by $i_{\gamma}$ the interior multiplication. The next algebraic
fact follows by direct computations
\begin{pro}\label{pro7}
For any 1-form $\gamma$ the identity \hspace{3mm}
$*(\Phi\wedge\gamma)=i_{\gamma}\Phi $ \hspace{3mm} holds.
\end{pro}
\begin{th}\label{thf}
Let $(M,g,\Phi)$ be a $Spin(7)$ manifold with fundamental 4-form
$\Phi$, $P$ be the corresponding 3-fold vector cross product and
$\nabla^g$ be the Levi-Civita connection of $g$. Then the
following formula holds for all vector fields $X,Y,Z,V,W$:
\begin{eqnarray}\label{n1}
\nonumber (\nabla^g_X\Phi)(Y,Z,V,W)&=&\\ \nonumber & &
\frac{1}{2}\left\{\delta \Phi(X,Y,P(Z,V,W)) - \delta \Phi(X,Z,P(Y,V,W))
\right\}\\  &+&\frac{1}{2}\left\{
\delta \Phi(X,V,P(Y,Z,W))- \delta \Phi(X,W,P(Y,Z,V))\right\}\\ \nonumber
&-&\frac{1}{12}\left\{*(\delta \Phi\wedge\Phi)(P(X,Y,P(Z,V,W)) -
*(\delta \Phi\wedge\Phi)(P(X,Z,P(Y,V,W))\right\}\\ \nonumber &+&
\frac{1}{12}\left\{ *(\delta \Phi\wedge\Phi)(P(X,V,P(Y,Z,W))-
*(\delta \Phi\wedge\Phi)(P(X,W,P(Y,Z,V))\right\}
\end{eqnarray}
\end{th}
{\bf Proof.} We have the general formulas (see eg \cite{KN})
\begin{eqnarray}\label{n2}
(\nabla^g_X\Phi)(Y,Z,V,W)&=& X\Phi(Y,Z,V,W)-\Phi(\nabla^g_XY,Z,V,W) -
\Phi(Y,\nabla^g_XZ,V,W)\\ \nonumber
&-& \Phi(Y,Z,\nabla^g_XV,W) - \Phi(Y,Z,V,\nabla^g_XW),
\end{eqnarray}
\begin{eqnarray}\label{n3}
2g(\nabla^g_XY,Z)&=& Xg(Y,Z)+Yg(X,Z)-Zg(X,Y)\\ \nonumber
&+& g([X,Y],Z)+g([Z,X],Y)-g([Y,Z],X).
\end{eqnarray}
Let $\{e_0,e_1,...,e_7\}$ be an orthonormal basis and the
fundamental form $\Phi$ be given by (\ref{1}). We substitute
(\ref{n3}) into (\ref{n2}). Using the expression (\ref{1}) and
keeping in mind Proposition~\ref{pro7}, we check that the right
hand side of the obtained equality coincides with the right hand
side of (\ref{n1}) by long but straightforward calculations
evaluating the both sides on the basis  $e_0,e_1,...,e_7$. \hfill
{\bf Q.E.D.}

\section{Proof of Theorem~\ref{th2} part i)}

Suppose that a connection $\nabla$ determined by \beq\label{n5}
g(\nabla_XY,Z)=g(\nabla^g_XY,Z)+\frac{1}{2}T(X,Y,Z), \eeq where
$T$ is a 3-form, satisfies $\nabla\Phi=0$. Then we have
\begin{eqnarray}\label{n6}
2(\nabla^g_X\Phi)(Y,Z,V,W)&=& \Phi(T(X,Y),Z,V,W) +
\Phi(Y,T(X,Z),V,W)\\ \nonumber &+& \Phi(Y,Z,T(X,V),W)+
\Phi(Y,Z,V,T(X,W))
\end{eqnarray}
and consequently \beq\label{n7} \delta\Phi = -*d*\Phi =
\sum_{i,j=0}^{7}\left((i_{e_j}i_{e_i}T)\wedge (i_{e_j}i_{e_i}\Phi)
\right) \eeq Evaluating (\ref{n7}) on the orthonormal basis and
using the expression of the fundamental 4-form (\ref{1}) with
respect to this basis we arrive to a linear system of maximal rank
of 56 linear equations with respect to 56 unknown variables
$T(e_i,e_j,e_k), i,j,k=0,..,7$ since $T$ is a 3-form. By the
symmetries of the fundamental 4-form this system is separated into
8 linear systems and each of them consists of 7 linear equations
with respect to 7 unknown variables. Solving each of these systems
explicitly and using the definition of the Lee form $\theta$ we
obtain (\ref{n4}).

For the converse, we define by (\ref{n5}) a connection $\nabla$
with totally skew symmetric torsion $T$ given by (\ref{n4}).
Clearly $\nabla g=0$. Substitute (\ref{n4}) into (\ref{n6}) and
using Theorem~\ref{thf} we get $\nabla\Phi=0$.

Let $\phi$ be the spinor corresponding to $\Phi$. Clearly $\phi$
is $\nabla$ parallel. The Clifford action $T\cdot\phi$ depends
only on the $\Lambda^3_8$-part of $T$. Using (\ref{n4}) and the
algebraic formulas $*(\gamma\wedge\Phi)\cdot\phi=
i_{\gamma}(\Phi)\cdot\phi = 7\gamma\cdot\phi$ we obtain
(\ref{b2}). This proves part i). Part ii) will be proved in the
next section. \hfill {\bf Q.E.D.}

Further, we shall call the connection determined by
Theorem~\ref{th2}  {\it the $Spin(7)$- connection} of a given
$Spin(7)$ manifold.
\begin{co}\label{co2}
The Lee 1-form of any  $Spin(7)$ structure and the projections
$\pi^3_8(d\Phi)$, $\pi^3_8(T)$ onto the space $\Lambda^3_8$ are
given by $ \theta = \frac{6}{7}*(\Phi\wedge T), \quad
\pi^3_8(d\Phi)=\theta\wedge\Phi, \quad
\pi^3_8(T)=-\frac{1}{6}*\left(\theta\wedge\Phi\right).$
\end{co}

Keeping in mind Proposition~\ref{pro1}, we get
\begin{co}
The torsion 3-form $T$ of the $Spin(7)$ connection $\nabla$
changes by a conformal transformation $(g_o=e^{2f}g,
\Phi_o=e^{4f}\Phi)$ of the $Spin(7)$ structure $(g,\Phi)$ by
$
T_o=e^{4f}\left(T - \frac{2}{3} *(df\wedge\Phi)\right).
$
\end{co}

\section{The Ricci tensor and the scalar curvature}

In this section we give formulas for the Ricci tensor and the
scalar curvature of the connection $\nabla$ on a $Spin(7)$
manifold and, consequently, formulas for the Ricci tensor and the
scalar curvature of the metric $g$ using the special properties of
the Clifford action on the $\nabla$-parallel spinor. We apply the
Schr\"odinger-Lichnerowicz formula for the Dirac operator of a
metric connection with totally skew-symmetric torsion proved in
\cite{FI} to the case of the unique $Spin(7)$-connection $\nabla$
on a $Spin(7)$ manifold $(M,g,\Phi)$. Finally, we prove the part
ii) of Theorem~\ref{th2}.

Let $D, Ric, Scal$ be the Dirac operator, the Ricci tensor and the
scalar curvature of the $Spin(7)$ connection defined as usually by
$D=\sum_{i=0}^7 e_i.\nabla_{e_i},\quad
Ric(X,Y)=\sum_{i=0}^7R(e_i,X,Y,e_i), \quad Scal =
\sum_{i=0}^7Ric(e_i,e_i).$ The relations between the Ricci tensor
$Ric^g$ and the scalar curvature $Scal^g$ of the metric are (see
\cite{IP,FI}) \beq\label{c5} Ric^g=Ric +\frac{1}{2}\delta T
+\frac{1}{4}(i_{(.)}T,i_{(.)}T), \quad Scal^g=Scal
+\frac{1}{4}||T||^2, \eeq where $(,)$ and $ ||.||^2$ denote the
inner product on tensors induced by $g$ and the corresponding
norm. In particular, $Ric$ is symmetric iff the torsion 3-form is
coclosed, $\delta T=0$.

Let $\sigma^T$ be the 4-form defined by $\sigma^T=\frac{1}{2}
\sum_{i=0}^7(i_{e_i}T)\wedge (i_{e_i}T)$. We take the following
result from \cite{FI}.

\begin{th}\label{tA}\cite{FI}
Let $\Psi$ be a parallel spinor with respect to a metric
connection $\nabla$ with totally skew- symmetric torsion $T$ on a
Riemannian spin manifold $M$. The following formulas hold
\begin{eqnarray}\label{c3}
& &\nonumber 3dT \cdot \Psi -2\sigma^T\cdot\Psi
+Scal\hspace{1mm}\Psi=0,\\ & &\frac{1}{2}i_XdT\cdot\Psi
+\nabla_XT\cdot\Psi -Ric(X)\cdot\Psi=0,\\ &
&D(T\cdot\Psi)=dT\cdot\Psi +\delta T\cdot\Psi
-2\sigma^T\cdot\Psi.\nonumber
\end{eqnarray}
If $M$ is compact, then for any spinor field $\psi$ the following
formula is true \beq\label{c4} \int_M||D\psi||^2\,dVol=
\int_M\left(||\nabla\psi||^2+
(dT\cdot\psi,\psi)+2(\sigma^T\cdot\psi,\psi) +
Scal||\psi||^2\right)\,dVol. \eeq In particular, if the
eigenvalues of the endomorphism $dT+2\sigma^T + Scal$ acting on
spinors are nonnegative, then every $\nabla$-harmonic spinor is
$\nabla$-parallel. If the eigenvalues are positive, then there are
no $\nabla$-parallel spinors.
\end{th}

We apply Theorem~\ref{tA} to the $\nabla$-parallel spinor $\phi$
corresponding to the fundamental 4-form $\Phi$ on a $Spin(7)$
manifold to get
\begin{pro}\label{thc1}
Let $(M,g,\Phi,\nabla)$ be an 8-dimensional $Spin(7)$ manifold
with the $Spin(7)$ connection $\nabla$ of torsion $T$. The Ricci
tensors $Ric,Ric^g$ are given by \beq\label{ric1}
Ric(X)=-\frac{1}{2} *\left(i_XdT\wedge\Phi\right)
-*\left(\nabla_XT\wedge\Phi\right), \eeq \beq\label{rg1}
Ric^g(X,Y)=\frac{1}{2} \left(i_XdT\wedge\Phi,*Y\right)
+\left(\nabla_XT\wedge\Phi,*Y\right) +\frac{1}{2}\delta T(X,Y)
+\frac{1}{4}(i_XT,i_YT)\quad . \eeq
\end{pro}

\subsection{Proof of Theorem~\ref{th2} ii)}

Let $\phi$ be the $\nabla$-parallel spinor corresponding to the
fundamental 4-form $\Phi$. Then the Riemannian Dirac operator
$D^g$ and the Levi-Civita connection $\nabla^g$ act on $\phi$ by
the rule \beq\label{dir1} \nabla^g_X\phi
=-\frac{1}{4}(i_XT)\cdot\phi, \quad D^g\phi = -
\frac{3}{4}T\cdot\phi=\frac{7}{8}\theta\cdot\phi, \eeq where we
used (\ref{b2}). We are going to apply the well known
Schr\"odinger-Lichnerowicz (S-L) formula \cite{Li,Sro}
\hspace{3mm} $(D^g)^2= \triangle^g  +\frac{1}{4}Scal^g, \quad
\triangle^g= -\sum\left( \nabla^g_{e_i}\nabla^g_{e_i}
-\nabla^g_{\nabla_{e_i}e_i}\right)$ to the  $\nabla$-parallel
spinor field $\phi$.

Using (\ref{dir1}) we calculate as a consequence that
\beq\label{d2} (D^g)^2\phi = \frac{7}{8}D^g(\theta\phi)=
\left(\frac{49}{64}||\theta||^2+
\frac{7}{8}\delta\theta\right)\cdot\phi+\frac{7}{8}d\theta\cdot\phi
+\frac{7}{16}(i_{\theta}T)\cdot\phi, \eeq where we used the
general identity $D^g\theta+\theta
D^g=d\theta+\delta\theta-2\nabla_{\theta}$.

We compute the Laplacian $\triangle^g$ in the general.
\begin{lem}\label{lq}
Let $\phi$ be a parallel spinor with respect to a metric
connection $\nabla$ with skew symmetric torsion $T$ on a
Riemannian manifold $(M,g)$. For the Riemannian Laplacian acting
on $\phi$ we have \beq\label{d3}
\triangle^g\phi=-\frac{1}{4}\delta T\cdot\phi
-\frac{1}{16}\left(2\sigma^T-\frac{1}{2}||T||^2
\right)\cdot\phi.\eeq
\end{lem}
{\bf Proof of Lemma~\ref{lq}.} We take a normal coordinate system
such that $(\nabla_{e_i}e_i)_p=0, p\in M$. We use (\ref{dir1}) to
get
$\triangle^g\phi=\frac{1}{4}\sum_i\left(\nabla_{e_i}i_{e_i}T)\cdot\phi
-\frac{1}{16}(i_{e_i}T)\cdot (i_{e_i}T)\cdot\phi\right).$ Applying
the properties of the Clifford multiplication we obtain (\ref{d3})
and Lemma~\ref{lq} is proved.

Further, substituting (\ref{d2}) and (\ref{d3}) into the S-L
formula, multiplying the obtained result by $\phi$ and taking the
real part, we arrive at \beq\label{d4}
\left(\frac{49}{64}||\theta||^2+\frac{7}{8}\delta\theta\right)||\phi||^2=
\left(\frac{1}{32}||T||^2+\frac{1}{4}Scal^g\right)||\phi||^2
-\frac{1}{8}(\sigma^T\cdot\phi,\phi). \eeq On the other hand,
using (\ref{b2}), we get
$D(T\cdot\phi)=-\frac{7}{6}D(\theta\cdot\phi)= -\frac{7}{6}
\left(d^{\nabla}\theta \cdot\phi + \delta\theta \cdot\phi\right),$
where $d^{\nabla}$ is the exterior derivative with respect to the
$Spin(7)$ connection $\nabla$. Now,  (\ref{c3}) gives \\
$-\frac{7}{6} \left(d^{\nabla}\theta \cdot\phi + \delta\theta\cdot
\phi\right) = dT\cdot\phi -2\sigma^T\cdot\phi +\delta T\cdot\phi$.
Multiplying the last equality by $\phi$ and taking the real part,
we obtain $-\frac{7}{6}\delta\theta ||\phi||^2=(dT\cdot\phi,\phi)
-(2\sigma^T\cdot\phi,\phi)$. Consequently,
 (\ref{c3}) and (\ref{c5}) imply
\beq\label{d5}
\left(-\frac{7}{2}\delta\theta-\frac{1}{4}||T||^2+Scal^g\right)||\phi||^2+
4(\sigma^T\cdot\phi,\phi)=0. \eeq Finally, we get (\ref{sc1}) from
(\ref{d4}) and (\ref{d5}). Thus, the proof of Theorem~\ref{th2} is
completed. \hfill {\bf Q.E.D.}
\begin{co}
On a balanced $Spin(7)$ manifold the Ricci tensor $Ric$ is
symmetric and the Riemannian Ricci tensor and scalar curvatures
are given by
\begin{eqnarray}\label{ric2}\nonumber
& &Ric(X,Y)=\frac{1}{2} \left(i_X(d\delta\Phi)\wedge\Phi,*Y\right), \quad
Scal=-\frac{1}{3}||\delta\Phi||^2;\\
& &Ric^g(X,Y)=\frac{1}{2} \left(i_X(d\delta\Phi)\wedge\Phi,*Y\right)
 +\frac{1}{4}(i_XT,i_YT), \quad Scal^g=-\frac{1}{12}||\delta\Phi||^2.
\end{eqnarray}

In particular the Riemannian scalar curvature on a balanced
$Spin(7)$-manifold is non-positive and vanishes identically if and
only if the $Spin(7)$-structure is co-closed, $\delta\Phi=0$  and
therefore parallel.

A balanced $Spin(7)$-manifold has harmonic fundamental form,
$d\delta\Phi=0$ or equivalently it has closed torsion 3-form,
$dT=0$ if and only if the $Spin(7)$-structure is co-closed,
$\delta\Phi=0$  and therefore parallel.
\end{co}
{\bf Proof.} In the case of a balanced structure, the torsion
3-form $T$ satisfies $T=-\delta\Phi$ by Theorem~\ref{th2}.
Clearly, $\delta T=0$ and $Ric$ is a symmetric tensor. The
Clifford multiplication of a 3-form by the spinor $\phi$ depends
only on its projection in the space $\Lambda^3_8$. The 3-form  $T$
belongs to $\Lambda^3_{48}$ by Corollary~\ref{co2} and hence,
$\nabla T$, as a 3-form, also belongs to $\Lambda^3_{48}$ since
the $Spin(7)$ connection preserves the fundamental 4-form and
therefore it preserves also the splitting $\Lambda^p_l$. Hence,
the Clifford action of $\nabla T$ on the special spinor $\phi$ is
trivial. The rest of the claim follows from Theorem~\ref{th2} and
Proposition~\ref{thc1}. \hfill {\bf Q.E.D.}

\section{Topology of compact $Spin(7)$ manifold}

In this section we apply our results to obtain information about
Betti numbers, $\hat A$-genus and the signature of certain classes
of $Spin(7)$ manifolds. We use essentially the solution of the
Yamabe conjecture \cite{RS} as well as the fundamental
Atiyah-Singer Index theorem \cite{AS} which gives a topological
formula for the index of any linear elliptic operator. On a
$Spin(7)$ manifold $M$ this reads as $ind$ D$ =\hat A(M)=ind$$
D^g$, where $\hat A(M)$ is a topological invariant called $\hat
A$-genus, $ind$ D$=dim Ker D_+-dim KerD_-,
D_{\pm}:\Gamma(S_{\pm})\rightarrow \Gamma(S_{\mp})$ are the Dirac
operators of a linear connection on $M$.

First, we notice that the expression of the $\hat A$-genus in
terms of Betti numbers proved by Joyce \cite{J1,J2} for a parallel
compact $Spin(7)$ manifold holds for any compact $Spin(7)$
manifold.
\begin{pro}\label{top}
On a compact $Spin(7)$ manifold $(M,g,\Phi)$ the $\hat A$-genus is
given by
\beq\label{gen} 24\hat A(M) = -1+b_1-b_2+b_3+(b_4)_+
-2(b_4)_-, \eeq where $b_i$ are the Betti numbers of $M$ and
$(b_4)_+$ (resp. $(b_4)_-)$ is the dimension of the space of
harmonic self-dual (resp. anti-self dual) 4-forms.
\end{pro}
{\bf Proof.} The proof goes as in \cite{J2} following the
reasoning of \cite{Sal1}. We recall the basic identities. The
formula for the signature $\tau(M)$ and the $\hat A$-genus in
terms of Pontrjagin classes are \cite{Hir} \beq\label{sign}
45\left((b_4)_+-(b_4)_-\right)=45\tau(M)=7p_2(M)-p^2_1(M), \quad
45.2^7\hat A(M)=7p_1^2(M)-4p_2(M). \eeq Combining (\ref{sign})
with (\ref{c1}) gives (\ref{gen}). \hfill {\bf Q.E.D.}

We state the main result of this section
\begin{th}\label{mt}
Let $M$ be a compact connected spin 8-manifold with a fixed
orientation. If it admits a strict locally conformally parallel
$Spin(7)$ structure $(g,\Phi)$, then $M$ admits a Riemannian
metric $g_Y$ with strictly positive constant scalar curvature,
$Scal^{g_Y}>0$.

Consequently, the following formulas hold
\item
i). $\hat A(M)=0$;
\item
ii).${\cal X}(M)=3\tau(M)$;
\item
iii). $b_2+2(b_4)_+-b_3-(b_4)_-\ge 0$ with equality iff $b_1=1$.

In particular, $M$ does not admit a metric with holonomy
$Hol(g)=Spin(7); SU(4); Sp(2)$;  $SU(2)\times SU(2)$.
\end{th}
{\bf Proof.} Let $\theta$ be the Lee form of $(g,\Phi)$. We need
the following algebraic lemma.
\begin{lem}\label{al}
On a $Spin(7)$ manifold the inequality \hspace{3mm}
$
||T||^2\ge \frac{7}{6}||\theta||^2 $ \hspace{3mm} holds. The
equality is attained if and only if the $Spin(7)$ structure is
locally conformally parallel.
\end{lem}
The proof of Lemma~\ref{al} follows from Theorem~\ref{th2} and the  equality
$0\leq ||T+\frac{1}{6}*(\theta\wedge\Phi)||^2=
||T||^2-\frac{7}{6}||\theta||^2$.

Lemma~\ref{al} gives $||T||^2=\frac{7}{6}||\theta||^2$ since the
structure is locally conformally parallel. Theorem~\ref{th2} leads
to the formula \beq\label{as} Scal^g =\frac{21}{36}||\theta||^2
+\frac{7}{2}\delta\theta \quad . \eeq According to the solution of
the Yamabe conjecture \cite{RS} there is a metric $g_Y=e^{2f}g$ in
the conformal class of $g$ with constant scalar curvature,
$Scal^{g_Y}= const.$ Consider the locally conformally parallel
$Spin(7)$ structure $(g_Y=e^{2f}g, \Phi_Y=e^{4f}\Phi)$.  Equality
(\ref{as}) is true also for the structure $(g_Y,\Phi_Y)$. An
integration of the last equality over a compact $M$ gives $$
Scal^{g_Y}.Vol_{g_Y}= \frac{21}{36}\int_M||\theta||^2\,dVol_{g_Y}
>0, $$ since the structure is strictly locally conformally parallel.
Then, by the Lichnerowicz vanishing theorem \cite{Li}, $ind
D^{g_Y}=0$ and $\hat A(M)=0$ by the index theorem. Condition ii)
follows exactly as in \cite{Sal1} from (\ref{sign}) and
(\ref{c1}). Statement iii) is a consequence of (\ref{gen}) and
Corollary~\ref{cof}. We derive the last assertion  by
contradiction with the already proved vanishing of the $\hat
A$-genus and the result of Joyce \cite{J1,J2} claiming
non-vanishing of the $\hat A$-genus for a Riemannian manifold with
Riemannian holonomy groups listed in the condition of the theorem.
\hfill {\bf Q.E.D.}

{\bf Remark} Information for the $\hat A$-genus on a compact
$Spin(7)$ manifold can be obtained if the eigenvalues of the
endomorphism $dT+2\sigma^T+Scal$ acting on spinors are known
according to Theorem~\ref{tA}. In particular, if the eigenvalues
are non-negative (they cannot be positive since there always
exists a parallel spinor), then the holonomy group $Hol(\nabla)$
will determine the $\hat A$ genus in the simply connected case
since the index of $D$ is given by the $\nabla$-parallel spinors.
For example, if $Hol(\nabla)=Spin(7); SU(4); Sp(2); SU(2)\times
SU(2)$, then $\hat A=1; 2; 3; 4$, respectively by  pure algebraic
arguments, namely by considering the fixed spinors by the action
of the holonomy representation of $\nabla$ on spinors.

\section{Solutions to the Killing spinor equations in dimension 8}

We consider the Killing spinor equations (\ref{ks1}) and
(\ref{ks2}) in dimension 8. The existence of a non-trivial
$\nabla$-parallel spinor is equivalent to  the existence of  a
$Spin(7)$ structure $(g,\Phi)$ \cite{LM}. Then
 the 3-form field strength $H=T$ is given by Theorem~\ref{th2}.
Involving the second Killing spinor equation (\ref{ks2}) we have
\begin{th}\label{th3}
In dimension 8 the following conditions are equivalent:

i) The Killing spinor equations (\ref{ks1}) and (\ref{ks2}) admit
solution with dilation $\Psi$;

ii) There exists a $Spin(7)$ structure $(g,\Phi)$ with closed Lee
form $\theta=-\frac{12}{7}d\Psi$ and therefore it is locally
conformal to a balanced $Spin(7)$ structure.

The 3-form field strength  $H=T$ and the Riemannian scalar
curvature $Scal^g$ are given by \beq\label{6} T=-\delta\Phi  +
2*(d\Psi\wedge\Phi), \eeq \beq\label{lsc} Scal^g =8||d\Psi||^2
-\frac{1}{12}||T||^2-6\triangle\Psi, \eeq where $\triangle
\Psi=\delta d\Psi$ is the Laplacian.

The solution is with constant dilation if and only if the
$Spin(7)$ structure is balanced.
\end{th}
{\bf Proof.} We apply Theorem~\ref{th2}. Let $\nabla$ be a
connection with torsion 3-form $T$. Let $\phi$ be an arbitrary
$\nabla$-parallel spinor field such that $(2d\Psi-T)\cdot\phi=0$.
The spinor field $\phi$ defines a Spin(7) structure $\Phi$ which
is $\nabla$-parallel. On the other hand, the connection preserving
$\Phi$ with torsion any 3-form is unique given by
Theorem~\ref{th2}. Comparing (\ref{b2}) with the second Killing
spinor equation (\ref{ks2}) we find $\frac{12}{7}d\Phi=-\theta$.
Inserting the last equality into (\ref{n4}) and (\ref{sc1}), we
get (\ref{6}) and (\ref{lsc})
 which completes the proof. \hfill {\bf Q.E.D.}

A similar formula as (\ref{6})  was derived in \cite{GKMW}
as a necessary condition.

Theorem~\ref{th3} allows us to obtain a lot of compact solutions
to the Killing spinor equations. If the dilation is a globally
defined function, then any solution is globally conformal
equivalent to a balanced $Spin(7)$ structure. For example, any
conformal transformation of a compact 8-dimensional manifold with
Riemannian holonomy group $Spin(7)$ constructed by Joyce
\cite{J1,J2} is a solution with a globally defined non-constant
dilation.

Summarizing we obtain
\begin{co}
Any solution $(M^8,g,\Phi)$ to the Killing spinor equations
(\ref{ks1}), (\ref{ks2}) in dimension 8 with non-constant globally
defined dilation function $\Psi$ comes from a solution with
constant dilation by a conformal transformation ie
$(g=e^{\frac{6}{7}\Psi}g_0,\Phi=e^{\frac{12}{7}\Phi}_0)$, where
$(g_0,\Phi_0)$ is a balanced $Spin(7)$ structure.
\end{co}

{\bf Note added to the proof.} It has been shown in \cite{Fr} that
if an $n$-ddimensional $G$-structure with structure group $G$
satisfying certain weak conditions admits a $G$-connection with
totally skew-symmetric torsion then the $G$-structure has to be a
Spin(7)-structure in dimension 8.

\vspace{5mm}
\noindent University of Sofia "St. Kl. Ohridski"   \\
Faculty of Mathematics and Informatics,\\
blvd. James Bourchier 5,\\
1164 Sofia, Bulgaria\\
\end{document}